\documentclass[11pt,a4paper]{article}

\usepackage{a4wide,amssymb,graphics,graphicx,textcomp}
\usepackage{enumerate,textcomp,multirow,amsmath}
\usepackage{algorithm}
\usepackage{algorithmic}
\usepackage{url}

\newtheorem{theorem}{Theorem}[section]
\newtheorem{lemma}[theorem]{Lemma}
\newtheorem{remark}[theorem]{Remark}

\newtheorem{corollary}[theorem]{Corollary}

\numberwithin{equation}{section}
\numberwithin{table}{section}
\numberwithin{figure}{section}
\newcommand{\proof}{\par\noindent{\bf Proof}. \ignorespaces}
\newcommand{\eproof}{\space
    {\ \vbox{\hrule\hbox{\vrule height1.3ex\hskip0.8ex\vrule}\hrule}}\par}

\newcommand {\mat}  [1] {\left[\begin{array}{#1}}
\newcommand {\rix}      {\end{array}\right]}

\catcode`@=11     
\font\tenex=cmex10 
\newdimen\p@renwd
\setbox0=\hbox{\tenex B} \p@renwd=\wd0 
\def\bmat#1{\begingroup \m@th
  \setbox\z@\vbox{\def\cr{\crcr\noalign{\kern2\p@\global\let\cr\endline}}%
    \ialign{$##$\hfil\kern2\p@\kern\p@renwd&\thinspace\hfil$##$\hfil
      &&\quad\hfil$##$\hfil\crcr
      \omit\strut\hfil\crcr\noalign{\kern-\baselineskip}%
      #1\crcr\omit\strut\cr}}%
  \setbox\tw@\vbox{\unvcopy\z@\global\setbox\@ne\lastbox}%
  \setbox\tw@\hbox{\unhbox\@ne\unskip\global\setbox\@ne\lastbox}%
  \setbox\tw@\hbox{$\kern\wd\@ne\kern-\p@renwd\left[\kern-\wd\@ne
    \global\setbox\@ne\vbox{\box\@ne\kern2\p@}%
    \vcenter{\kern-\ht\@ne\unvbox\z@\kern-\baselineskip}\,\right]$}%
  \null\;\vbox{\kern\ht\@ne\box\tw@}\endgroup}
\catcode`@=12    

\newcommand{\argmin}{\textrm{argmin}}

\newcommand{\backmatter}{\appendix
\def\chaptermark##1{\markboth{%
\ifnum  \c@secnumdepth > \m@ne  \@chapapp\ \thechapter:  \fi  ##1}{%
\ifnum  \c@secnumdepth > \m@ne  \@chapapp\ \thechapter:  \fi  ##1}}%
\def\sectionmark##1{\relax}}

\makeatletter
\newcommand*{\rom}[1]{\expandafter\@slowromancap\romannumeral #1@}
\makeatother

\def\rank{\mathop{\mathrm{rank}}}

\def\diag{\mathop{\mathrm{diag}}}

\newif\ifMDlatex
\catcode`@=11     
\def\MD@us#1{\csname#1style\endcsname}
\def\MD@uf#1{\csname#1font\endcsname}
\def\MD@t{text}
\def\MD@s{script}
\def\MD@ss{scriptscript}
\newdimen\MD@unit
\MD@unit\p@
\def\MD@changestyle#1{
  \relax\MD@unit0.1\fontdimen6\MD@uf{#1}0
  \everymath\expandafter{\the\everymath\MD@us{#1}}
}
\def\MD@dot{$\m@th\ldotp$}
\def\MD@palette#1{\mathchoice{#1\MD@t}{#1\MD@t}{#1\MD@s}{#1\MD@ss}}
\def\MD@ddots#1{{\MD@changestyle{#1}%
  \mkern1mu\raise7\MD@unit\vbox{\kern7\MD@unit\hbox{\MD@dot}}%
  \mkern2mu\raise4\MD@unit\hbox{\MD@dot}%
  \mkern2mu\raise \MD@unit\hbox{\MD@dot}\mkern1mu}}%
\def\MD@iddots#1{{\MD@changestyle{#1}%
  \mkern1mu\raise \MD@unit\hbox{\MD@dot}%
  \mkern2mu\raise4\MD@unit\hbox{\MD@dot}%
  \mkern2mu\raise7\MD@unit\vbox{\kern7\MD@unit\hbox{\MD@dot}}}}%
\def\MD@vdots#1{\vbox{\MD@changestyle{#1}%
    \baselineskip4\MD@unit\lineskiplimit\z@
    \kern6\MD@unit\hbox{\MD@dot}\hbox{\MD@dot}\hbox{\MD@dot}}}%
\ifMDlatex
  \DeclareRobustCommand\ddots{\mathinner{\MD@palette\MD@ddots}}%
  \DeclareRobustCommand\iddots{\mathinner{\MD@palette\MD@iddots}}%
  \DeclareRobustCommand\vdots{\mathinner{\MD@palette\MD@vdots}}%
\else
  \def\ddots{\mathinner{\MD@palette\MD@ddots}}%
  \def\iddots{\mathinner{\MD@palette\MD@iddots}}%
  \def\vdots{\mathinner{\MD@palette\MD@vdots}}%
\fi
\catcode`@=12    

\newcommand {\comment}[1]{} 
\newcommand{\C}{{\mathbb C}}
\newcommand{\R}{{\mathbb R}}

\begin{document}
\title{\hspace{1.3cm} A semi-analytical approach for the positive semidefinite Procrustes problem}
\author{  Nicolas Gillis$^*$ \qquad Punit Sharma\thanks{The authors acknowledge the support of the ERC (starting grant n$^\text{o}$ 679515).
NG also acknowledges the support of the F.R.S.-FNRS (incentive grant for scientific research n$^\text{o}$ F.4501.16). Email: \{nicolas.gillis, punit.sharma\}@umons.ac.be.} \\ 
Department of Mathematics and Operational Research \\
Facult\'e Polytechnique, Universit\'e de Mons \\
Rue de Houdain 9, 7000 Mons, Belgium
}

\date{}

\maketitle

\begin{abstract}
The positive semidefinite Procrustes (PSDP) problem is the following: given
rectangular matrices $X$ and $B$, find the symmetric positive semidefinite
matrix $A$ that minimizes the Frobenius norm of $AX-B$.
No general procedure is known that gives an exact solution.
In this paper, we present a semi-analytical approach to solve the PSDP problem.
First, we characterize a family of positive semidefinite matrices that either solve the PSDP problem
when the infimum is attained or give arbitrary accurate approximations to the infimum
when it is not attained.
This characterization requires the unique optimal solution of a smaller PSDP problem where $B$ is square and $X$ is diagonal
with positive diagonal elements.
Second, we propose a very efficient strategy to solve the PSDP problem,
combining the semi-analytical approach, a new initialization strategy and the fast gradient method.
We illustrate the effectiveness of the new approach, which is guaranteed to converge linearly,
compared to state-of-the-art methods.
\end{abstract}

\section{Introduction}
Given $n \times m$ matrices $X$ and $B \in \R^{n,m}$, the corresponding
positive semidefinite Procrustes (PSDP) problem is defined by
\begin{equation} \label{psdproc}
\inf_{ A \in \mathcal{S}_{\succeq}^n } \|AX-B\|_F^2, \quad\quad \text{where}~X,B \in \R^{n,m}, \tag{$\mathcal{P}$}
\end{equation}
where ${\|\cdot\|}_F$ denotes the Frobenius norm of a matrix and ${\mathbb R}^{n,r}$
denotes the set of $n \times r$ real matrices with the special case
${\mathbb R}^n={\mathbb R}^{n,1}$. The set of symmetric positive semidefinite matrices of size $n$ is
denoted by $\mathcal S_{\succeq}^n$.

This problem occurs for example in the field of structure analysis~\cite{Bro68} and in signal
processing~\cite{SufH93}.
For an elastic structure, each column of $X$ consists of generalized forces while each column of $B$ consists of the corresponding displacements, for a set of $m$-measurements. From this data, it is possible to recover the so-called compliance matrix $A$ that relates these column vectors by $AX=B$ and that must be symmetric positive definite. Such a compliance matrix may not exist for
the available measurements and it is therefore desirable to find the matrix $A$ that solves~\eqref{psdproc} instead~\cite{Woo96}.
Solutions to~\eqref{psdproc} can also be used when 
looking for the nearest stable matrix to an unstable one
using a block-coordinate descent method~\cite{GS16}. This is what initially motivated us to study this problem.

In the simplest case when $X$ is the identity matrix, the nearest positive semidefinite matrix to $B$
in the Frobenius norm is given by $(C+H)/2$, where $H$ is the symmetric polar factor of the matrix \mbox{$C=(B+B^T)/2$} \cite{Hig88b}. Equivalently, the projection $\mathcal P_{ \succeq }(B)$ of $B$ onto the cone of semidefinite matrices is given by
\begin{equation} \label{projpsd}
\mathcal P_{ \succeq}(B)=U\left(\max{(\Gamma,0)}\right)U^T,
\end{equation}
where $U \Gamma U^T$ is an eigenvalue decomposition of the symmetric matrix $\frac{B+B^T}{2}$.

The problem of finding the nearest Hermitian positive semidefinite matrix with a Toeplitz structure is studied in~\cite{SufH93}.
If the feasible set in~\eqref{psdproc} is chosen to be the set of orthogonal matrices,
then the problem is called the orthogonal Procrustes problem which arises in many applications such as computer vision, factor analysis, multidimensional scaling, and manifold optimization; see~\cite{GowD04, Gre52, Sch66, AbsMal2012} and references therein.
On the other hand, the symmetric Procrustes problem, where the feasible set in~\eqref{psdproc}
is the set of symmetric matrices arise in applications such as determination of
space craft altitudes and the stiffness matrix of an elastic structure~\cite{Hig88c,Bro68,Lar66}.

In a recent work, Alam and Adhikari~\cite{BibA16} have characterized and determined all solutions
of the structured Procrustes problem analytically, where the
feasible set in~\eqref{psdproc} is either a Jordan algebra or a Lie algebra associated with
an appropriate scalar product on $\R^n$ or $\C^n$.

The theoretical and computational aspects of PSDP problem have been
extensively studied in the literature.
It was first introduced and studied by Allwright~\cite{All88}, and later in detailed
by Woodgate~\cite{Woo87, Woo96,Woo98}, see also~\cite{SufH93, AE97, krislock2004local, krislock2003numerical, DenHZ03, Lia99} for more references.
A necessary and sufficient condition was provided for the existence of a solution for the PSDP problem by
Allwright and Woodgate in~\cite{AllW90}. An expression for the solution of the PSDP problem
was determined in~\cite{Lia99,Woo96} for some special cases. However,
we note that no general procedure is known for solving~\eqref{psdproc} analytically. Many algorithmic
solutions have been proposed in the literature, see for example~\cite{All88, Woo96, Woo98, AE97, krislock2003numerical}.
Other closely related problems are the computation of the positive semidefinite least square solutions of the matrix equation $A^TXA=B$~\cite{HuaDL96} and the equation $AXB = C$ in variable $X$~\cite{LiaB03}, and complex PSDP problems~\cite{KisH07}.

\paragraph{Notation}In the following, ${\|\cdot\|}_2$ denotes the spectral norm of a vector or a matrix.
For $A=A^T \in \R^{n,n}$, we denote $A\succ 0$ and $A\succeq 0$ if $A$ is symmetric positive definite
or symmetric positive semidefinite, respectively. By $\sigma_i(X)$ we denote the $i$th largest singular value of $X$,
 and $A^\dag$ denotes the Moore-Penrose pseudoinverse of the matrix $A$.

We use the acronym SVD for the singular value decomposition.

\subsection{Contributions and outline of the paper} \label{sec:contri}

In section 2, we derive some preliminary lemmas.
Here, we also provide a simpler alternative proof of the fact that the infimum in~\eqref{psdproc}
is always uniquely attained when $X$ has rank $n$.

In section 3, we present a semi-analytical solution to solve the PSDP  problem~\eqref{psdproc}, see Theorem~\ref{thm:general_analytic}.
 We reduce the original problem~\eqref{psdproc} into a smaller diagonal PSDP problem
 that always has a unique solution. We note that a similar reduction process
has been considered in~\cite[Theorem 2.1]{Lia99}.  As in~\cite[Theorem 2.1]{Lia99}, Theorem~\ref{thm:general_analytic}
 gives a necessary and sufficient condition for the infimum to be attained in $(\mathcal P)$, and characterizes the family of positive
semidefinite matrices that attains the infimum assuming the solution for the smaller diagonal PSDP problem is known.
The contribution of this section is to provide a more comprehensive description of the solutions of the PSDP problem:
\begin{itemize}

\item In addition to the conclusion of~\cite[Theorem 2.1]{Lia99},
Theorem~\ref{thm:general_analytic} completely characterizes the case \textit{when the infimum is not attained}
poviding the value of the infimum in $(\mathcal P)$ and deriving a family of arbitrary accurate approximations to the
infimum in~\eqref{psdproc}.

\item Because Theorem~\ref{thm:general_analytic} characterizes the set of optimal solutions of~\eqref{psdproc}, it can be used to obtain
 explicit formula for solutions to the PSDP problem, or arbitrary close approximations when the infimum is not attained,
with extremal properties of minimal rank, minimal Frobenius norm, or minimal spectral norm (Corollary~\ref{cor:corollary1}).

\item When $\text{rank}(X)=n$ and $(BX^T+XB^T) \preceq 0$, \cite[Theorem~2.5]{Woo96} showed that
the unique optimal solution in~$(\mathcal P)$ is the zero matrix.
When $\text{rank}(X) < n$ and  $U_1^T(BX^T+XB^T)U_1 \preceq 0$ with $U_1 \in \mathbb{R}^{n \times r}$ an orthogonal basis of the column space of $X$, the infimum is not attained in $(\mathcal P)$.
However, we explicitly provide the value of the infimum in $(\mathcal P)$ without computing the solution for the smaller diagonal PSDP problem,
and obtain analytically a family of arbitrary close approximations (Theorem~\ref{thm:particular}).

\item In the particular case $\text{rank}(X)=1$, 
we give a complete analytic characterization
for the optimal solutions in $(\mathcal P)$, see Theorem~\ref{thm:analytic_rank1}.

\end{itemize}

In section 4, we describe the fast gradient method (FGM) applied to~\eqref{psdproc}. FGM is an optimal first-order method and, in the strongly convex case, is guaranteed to converge linearly with optimal rate among first-order methods.
Therefore, combining the semi-analytical approach with FGM allows us to guarantee linear convergence of the objective function value with rate $(1-1/\kappa)$, where $\kappa = \frac{\sigma_1(X)}{\sigma_r(X)} \geq 1$ and $r = \rank(X)$.
Note that Andersson and Elfving~\cite{AE97} had already introduced the reduction of~\eqref{psdproc} but did not use it in their algorithms (also, they did not provide an explicit characterization of the solutions of~\eqref{psdproc} based on the solution of the subproblem, and their first-order algorithm is not optimal).
Moreover, to deal with ill-conditioned cases effectively (when $\kappa$ is large), we propose a very effective initialization strategy for~\eqref{psdproc} based on a recursive decomposition of the ill-conditioned problem into well-conditioned subproblems.

Finally in section 5, we present some numerical experiments illustrating the performance of the new proposed approach compared to state-of-the-art algorithms.

\section{Preliminary results}

If $X$ is of full rank, then the existence of a unique solution to the problem $(\mathcal P)$
is guaranteed by~\cite[Theorem~2.2]{Woo96}. A simplified proof of this fact was presented
in~\cite[Theorem 2.4.5]{krislock2003numerical} using a theorem of Weierstrass (``Suppose that the set $D$ is nonempty and closed, and that all the sub-level sets of the continuous function $f : D \rightarrow \mathbb{R}$ are bounded. Then $f$ has a global minimizer.'').
However, for the sake of completeness we restate the result and provide an alternative proof that uses the extreme value theorem (``a continuous function attains its infimum on a compact set'') and strong convexity, instead of Weierstrass' theorem.

\begin{lemma}\label{lem:psd_fullrank}
Let $X,B \in \R^{n,m}$. If $X$ has rank $n$, then the infimum of~\eqref{psdproc} is attained
for a unique $A\in \mathcal S_{\succeq}^{n}$.
\end{lemma}
\proof Let us define the continuous function $f : \mathcal S_{\succeq}^n \rightarrow \R $
such that $ A \rightarrow f(A)={\|AX-B\|}_F^2$.
Then $f$ is a strongly convex function as its  Hessian is positive definite: in fact, since $X$ has rank $n$
we have $\frac{d^2f}{dA^2}=\left(I_n \otimes XX^T \right) \succ 0$, where $I_n$ denotes the identity matrix of size $n$ and
$\otimes$ denotes the Kronecker product. Each diagonal block $XX^T$ in the Hessian corresponds to a row of $A$ since
${\|AX-B\|}_F^2 = \sum_{i=1}^n {\|A(i,:)X-B(i,:)\|}_F^2$.
Note that $\mathcal S_{\succeq}^n$ is a closed convex cone in the set of real symmetric matrices of size $n$, but not bounded.

However, the set
\[
\mathcal B_{\succeq}^n:=\left\{A\in \mathcal S_{\succeq}^n\,|~{\|AX-B\|}_F \leq {\|0 - B\|}_F = {\|B\|}_F\right \}
\]
is closed and bounded in $\mathcal S_{\succeq}^n$. In fact, the boundedness of the set
$\mathcal B_{\succeq}^n$ follows by \cite[Lemma~2.4.4]{krislock2003numerical}, because $A \in \mathcal B_{\succeq}^n$ implies that
\[
2 {\|B\|}_F \geq {\|AX\|}_F  \geq \frac{\|A\|_F}{\sqrt{n}}\, \sigma_n(X).
\]

Also, the zero matrix, that is, the matrix with all its entries
equal to zero, is an element of $\mathcal S_{\succeq}^n$, replacing the feasible set in~\eqref{psdproc} by $\mathcal B_{\succeq}^n$ does not change the infimum. Therefore, since $f$ is continuous, 
its infimum is attained at some point $\hat A \in \mathcal B_{\succeq}^n$. Further the uniqueness of such a $\hat A$ follows from the strong convexity of $f$.
\eproof

Following the strategy used in~\cite{BibA16} to derive analytic solutions of the structured Procrustes problem, we reduce
the original problem~\eqref{psdproc} to a smaller problem whose size equals the rank of $X$. The
following lemma that gives an equivalent characterization for a positive semidefinite
matrix will be repeatedly used in doing so.
\begin{lemma}[\cite{Alb69}] \label{lem:psd_character}
Let the integer $s$ be such that $0<s<n$, and $R=R^T \in \R^{n,n}$ be partitioned as
$R=\mat{cc}B & C^T\\C & D\rix$ with $B\in \R^{s,s}$, $C\in \R^{n-s,s}$ and $D \in \R^{n-s,n-s}$. Then $R \succeq 0$ if and only if
\begin{enumerate}
\item $B \succeq 0$,
\item $\operatorname{ker}(B) \subseteq \operatorname{ker}(C)$, and
\item $D-CB^{\dagger}C^T \succeq 0$, where $B^{\dagger}$ denotes the Moore-Penrose pseudoinverse of $B$.
\end{enumerate}
\end{lemma}

In the next section we use the fact that the trace of product of two positive semidefinite
matrices is nonnegative. The following elementary lemma gives more than that.

\begin{lemma}\label{lem:product_psd}
Let $P,Q \in \mathcal S_{\succeq}^n$. Then all eigenvalues of $PQ$ are nonnegative.
\end{lemma}
\proof Let $L \in \mathcal S_{\succeq}^n$ be such that $P=L^{\frac{1}{2}}L^{\frac{1}{2}}$ then
 $PQ=L^{\frac{1}{2}}(L^{\frac{1}{2}} Q)$ and $L^{\frac{1}{2}} Q L^{\frac{1}{2}}$ have the same nonzero
eigenvalues. But then $L^{\frac{1}{2}} Q L^{\frac{1}{2}} \in \mathcal S_{\succeq}^n$ as $Q \in \mathcal S_{\succeq}^n$
implies that all eigenvalues of the matrix $L^{\frac{1}{2}}(L^{\frac{1}{2}} Q)=PQ$ are nonnegative.
\eproof
We close the section with a result that will be used in obtaining solutions of
the PSDP problem~\eqref{psdproc} with the extremum properties of minimal rank, minimal Frobenius norm,
or minimal spectral norm, and investigate their uniqueness.

\begin{lemma}\label{lem:min_rank_F_2}
Let the integer $s$  be such that $0<s<n$. Let $B \in \mathcal S_{\succeq}^s$ and $C \in \R^{n-s,s}$ be such that
$\operatorname{ker}(B) \subseteq \operatorname{ker}(C)$. Define
$\mathcal D:=\{K \in \R^{n-s,n-s}\,:~K-CB^{\dag}C^T \succeq 0\}$ and  define $f:\mathcal D \rightarrow \R^{n,n}$ by
\[
f(K):=\mat{cc}B & C^T\\C & K\rix.
\]
Then the matrix $\hat K=CB^{\dag}C^T \in \mathcal D$ is a solution of  the  minimal rank
problem $\min_{K \in\mathcal D} \operatorname {rank}(f(K))$, the minimal Frobenius norm problem
$\min_{K \in\mathcal D} {\|f(K)\|}_F$, and
the minimal spectral  norm problem $\min_{K \in\mathcal D} {\|f(K)\|}_2$.
Moreover, for the minimal Frobenius norm and minimal rank problems, it is the unique solution.
\end{lemma}
\proof
Let $K \in \mathcal D$. Using Schur compliment of $B$ in $f(K)$, we have
\[
\operatorname{rank}(f(K))=\operatorname{rank}\left(\mat{cc}B&0\\0&K-CB^\dag C^T \rix\right)
\geq \operatorname{rank}(B).
\]
This implies that
\[
\min_{K \in \mathcal D}\operatorname{rank}(f(K))\geq \operatorname{rank}(B),
\]
and the minimum is uniquely attained when $K=CB^\dag C^T$.
For the minimal norm problems, observe that
\begin{eqnarray*}
\inf_{K\in \mathcal D}{\|f(K)\|}_{G}=
\inf_{K-CB^\dag C^T \succeq 0}
{\left\|\mat{cc} B& C^T\\ C & K\rix\right\|}_{G}
=\inf_{\Delta \succeq 0}
{\left\|\mat{cc} B & C^T\\ C & \Delta+CB^\dag C^T\rix
\right\|}_{G}, \qquad G = F \text{ or } 2.
\end{eqnarray*}
For the Frobenius norm, we have
\begin{eqnarray*}
\inf_{K\in \mathcal D}{\|f(K)\|}_{F}^2=\inf_{\Delta \succeq 0}
{\left\|\mat{cc}B & C^T\\ C & \Delta+CB^\dag C^T\rix
\right\|}_F^2
&=&
{\|B\|}_F^2+2{\|C\|}_F^2 + \inf_{\Delta \succeq 0}{\|\Delta+C{B}^{\dag}{C}^T\|}_F^2\\
&=& {\|B\|}_F^2+2{\|C\|}_F^2 + {\|C{B}^{\dag}{C}^T\|}_F^2,
\end{eqnarray*}
where the last equality holds as $\inf_{\Delta \succeq 0}{\|\Delta+C{B}^{\dag}{C}^T\|}_F$
is attained uniquely at $\Delta =0$ because $C{B}^{\dag}{C}^T \succeq 0$ as $B \succeq 0$.
Similarly, for the 2-norm,
\begin{eqnarray*}
\inf_{K\in \mathcal D}{\|f(K)\|}_{2}=\inf_{\Delta \succeq 0}
{\left\|\mat{cc}B & C^T\\ C & \Delta+CB^\dag C^T\rix
\right\|}_2 &=&
\inf_{\Delta \succeq 0}
{\bigg\|{\underbrace{\mat{cc}B & C^T\\ C &C{B}^{\dag}{C}^T\rix}_{R}}
+{\underbrace{\mat{cc}0 & 0\\ 0 & \Delta\rix}_{\Delta_R}} \bigg\|_2}\\
&=&
 \inf_{\Delta \succeq 0} \sup_{x \in \R^n \setminus \{0\}}\frac{x^*(R+\Delta_R)x}{x^*x} \geq
\sup_{x \in \R^n\setminus \{0\}}\frac{x^*Rx}{x^*x},
\end{eqnarray*}

where the last inequality is due to the fact that $x^*\Delta_Rx \geq 0$ for all $x \in \R^n$ as $\Delta_R \succeq 0$. Therefore the infimum is attained when $\Delta =0$.
\eproof

\section{Semi-analytical solutions for the PSDP  problem} \label{semianalsol}

The aim of this section is to provide comprehensive information about the solution(s) of the PSDP problem. 
As mentioned earlier,
we note that some of the results obtained here have partly appeared in~\cite{Lia99,Woo96} (in particular when the infimum is attained); see Section~\ref{sec:contri}.
Therefore,
we will relate the results in this section with those appeared in the literature by mentioning the
new contribution.

In the following, we present a semi-analytic solution for problem $(\mathcal P)$.
It is semi analytic in the sense that we reduce the original problem $(\mathcal P)$ into a
smaller problem that always has a unique solution. Then assuming the solution
for the subproblem is known, we characterize a family of positive semidefinite matrices that either
solve the problem~\eqref{psdproc} analytically or give arbitrary accurate approximations to the
infimum in~\eqref{psdproc}.

\begin{theorem}\label{thm:general_analytic}
Let $X,B \in \R^{n,m}$, and let $r = \operatorname{rank}(X)$. Let also $X=U\Sigma V^T$ be a singular value decomposition of $X$, where $U=\mat{cc}U_1 & U_2\rix \in \R^{n,n}$
with $U_1 \in \R^{n,r}$, $V=\mat{cc} V_1 & V_2\rix \in \R^{m,m}$ with $V_1 \in \R^{m,r}$,
and $\Sigma=\mat{cc} \Sigma_1 & 0\\0 & 0\rix \in \R^{n,m}$ with
$\Sigma_1\in \R^{r,r}$. Then
\begin{equation}\label{eq:psd_general_1}
\inf_{A\in \mathcal S_{\succeq}^n}{\|AX-B\|}_F^2=
\min_{A_{11}\in \mathcal S_{\succeq}^r}{\|A_{11}\Sigma_{1}-U_1^TBV_1\|}_F^2+{\|BV_2\|}_F^2.
\end{equation}
Further, let $\hat A_{11} \in \mathcal S_{\succeq}^r$ be such that
\begin{equation}\label{eq:subprob_A11}
\hat A_{11}=\operatorname{argmin}_{A_{11}\in \mathcal S_{\succeq}^r}{\|A_{11}\Sigma_{1}-U_1^TBV_1\|}_F^2.
\end{equation}
The following holds.
\begin{enumerate}
\item If $\operatorname{ker}(\hat A_{11}) \subseteq \operatorname{ker}(U_2^TBV_1\Sigma_1^{-1})$,
then $A_{opt}$ attains the infimum in~\eqref{eq:psd_general_1} if and only if
\begin{equation}\label{eq:psd_general_5}
A_{opt}:=U_1 \hat A_{11} U_1^T + U_2(U_2^TBV_1\Sigma_1^{-1})U_1^T+ U_1 (U_2^TBV_1\Sigma_1^{-1})^T U_2^T
+ U_2 K U_2^T,
\end{equation}
where $K\in \R^{n-r,n-r}$ is such that
$K-(U_2^TBV_1\Sigma_1^{-1})\hat A_{11}^{\dag}(U_2^TBV_1\Sigma_1^{-1})^T \succeq 0$.
\item Otherwise, the infimum in~\eqref{eq:psd_general_1} is not attained.
Let $\operatorname{rank}(\hat A_{11})= s < r$ and let $\epsilon > 0$ be sufficiently small.
Let $\hat A_{11}=\mat{cc}\hat U_1 & \hat U_2\rix \mat{cc}\hat \Sigma_1 & 0\\ 0 & 0 \rix \mat{c}\hat U_1^T\\ \hat U_2^T \rix$
be a singular value decomposition of $\hat A_{11}$, where $\hat U_1 \in \R^{r,s}$ and $\hat \Sigma_1 \in \R^{s,s}$.
Define
\begin{equation}\label{eq:psd_gen_6a}
\hat A_{11}^\epsilon:=\mat{cc}\hat U_1 & \hat U_2\rix \mat{cc}\hat \Sigma_1 & 0\\ 0& \Upsilon
\rix \mat{c} \hat U_1^T\\ \hat U_2^T \rix,
\end{equation}
where $\Upsilon \in \R^{r-s,r-s}$ is a diagonal matrix with diagonal entries each equal to
$\frac{\epsilon}{\beta}$, where
%
\begin{equation*}
\beta= \left\{\begin{array}{ll}
4\sqrt{(r-s)}{\|\Sigma_1\|}_F{\|\hat A_{11}\Sigma_1-U_1^TBV_1\|}_F & \mbox{ if }\,
{\|\hat A_{11}\Sigma_1-U_1^TBV_1\|}_F \neq 0, \\
4\sqrt{(r-s)}{\|\Sigma_1\|}_F & \mbox{ otherwise. } \end{array}\right.
\end{equation*}
Define
\begin{equation}\label{eq:psd_general_6}
A_{\epsilon}:=U_1 \hat A_{11}^\epsilon U_1^T + U_2(U_2^TBV_1\Sigma_1^{-1})U_1^T+ U_1 (U_2^TBV_1\Sigma_1^{-1})^T U_2^T
+ U_2 K_\epsilon U_2^T,
\end{equation}
where $K_\epsilon \in \R^{n-r,n-r}$ is such that
$K_\epsilon-(U_2^TBV_1\Sigma_1^{-1}) {(\hat A_{11}^\epsilon)}^{-1}(U_2^TBV_1\Sigma_1^{-1})^T \succeq 0$.
Then $A_\epsilon \in \mathcal S_{\succeq}^n$ and
\[
{\|A_\epsilon X-B\|}_F^2 < \inf_{A\in \mathcal S_{\succeq}^n}{\|AX-B\|}_F^2 + \epsilon.
\]
\end{enumerate}
\end{theorem}
\proof Let $A\in \R^{n,n}$ and set
\begin{equation}\label{eq:proof_general_1}
\hat A:=U^TAU=\mat{cc}A_{11} & A_{21}^T\\A_{21} & A_{22}\rix,
\end{equation}
where $A_{11}\in \R^{r,r}$, $A_{21}\in \R^{n-r,r}$ and $A_{22}\in \R^{n-r,n-r}$. Then
$A \succeq 0$ if and only if $\hat A \succeq 0$.
By Lemma~\ref{lem:psd_character}, $\hat A \succeq 0$ if and only if $A_{11}\succeq 0$,
$\operatorname{ker}(A_{11})\subseteq \operatorname{ker}(A_{21})$ and $A_{22}-A_{21}A_{11}^{\dag}A_{21}^T \succeq 0$.
Thus we have
\begin{eqnarray}\label{eq:proof_general_2}
{\|AX-B\|}_F^2&=&{\|U^TAU U^TX-U^TB\|}_F^2={\|\hat AU^TX-U^TB\|}_F^2 \nonumber\\
&=& {\left\|\mat{cc}A_{11} & A_{21}^T\\A_{21} & A_{22}\rix\mat{c}U_1^T\\U_2^T\rix X-
\mat{c}U_1^T\\U_2^T\rix B\right\|}_F^2 \nonumber \\
&=& {\left\|\mat{cc}A_{11} & A_{21}^T\\A_{21} & A_{22}\rix\mat{c}U_1^T X\\0\rix -
\mat{c}U_1^TB\\U_2^TB\rix \right\|}_F^2 \nonumber \\
&=& {\|A_{11}U_1^TX-U_1^TB\|}_F^2 +{\|A_{21}U_1^TX - U_2^T B\|}_F^2 \nonumber \\
&=& {\left \|A_{11}U_1^T\mat{cc}U_1 & U_2\rix\mat{cc}\Sigma_1 & 0 \\ 0& 0\rix-
U_1^TB\mat{cc}V_1 & V_2\rix\right\|}_F^2 +\nonumber \\
&& {\left\|A_{21}U_1^T\mat{cc}U_1 & U_2\rix\mat{cc}\Sigma_1 & 0 \\ 0& 0\rix-U_2^TB\mat{cc}V_1 & V_2\rix\right\|}_F^2 \nonumber \\
&=& {\left\|\mat{cc}A_{11}\Sigma_1-U_1^TBV_1 & -U_1^TBV_2\rix\right\|}_F^2 +
{\left\|\mat{cc}A_{21}\Sigma_1-U_2^TBV_1 & -U_2^TBV_2\rix\right\|}_F^2 \nonumber \\
&=& {\|A_{11}\Sigma_1-U_1^TBV_1\|}_F^2+{\|U_1^TBV_2\|}_F^2+{\|A_{21}\Sigma_1-U_2^TBV_1\|}_F^2+
{\|U_2^TBV_2\|}_F^2 \nonumber \\
&=& {\|A_{11}\Sigma_1-U_1^TBV_1\|}_F^2+ {\|A_{21}\Sigma_1-U_2^TBV_1\|}_F^2+{\|BV_2\|}_F^2,
\end{eqnarray}
where in the last equality we used ${\|U_1^TBV_2\|}_F^2+{\|U_2^TBV_2\|}_F^2={\|BV_2\|}_F^2$ since the
Frobenius norm is unitarily invariant. Thus taking the infimum in~\eqref{eq:proof_general_2} over $\mathcal S_{\succeq}^n$,
we obtain
\begin{align}
\inf_{A\in \mathcal S_{\succeq}^n} & \quad {\|AX-B\|}_F^2 \nonumber \\
& =  \inf_{A_{11}\succeq 0,\,
\operatorname{ker}(A_{11})\subseteq \operatorname{ker}(A_{21}),\,A_{22}-A_{21}A_{11}^{\dag}A_{21}^T \succeq 0}
{\|U\hat A U^T X-B\|}_F^2 \nonumber \\
& =  \inf_{A_{11}\succeq 0,\,\operatorname{ker}(A_{11})\subseteq \operatorname{ker}(A_{21})}
{\|A_{11}\Sigma_1-U_1^TBV_1\|}_F^2+ {\|A_{21}\Sigma_1-U_2^TBV_1\|}_F^2+{\|BV_2\|}_F^2, \label{eq:proof_general_3}
\end{align}
which implies that the infimum does not depend on the $A_{22}$ block of matrix $\hat A$ in~\eqref{eq:proof_general_1}.
Again from~\eqref{eq:proof_general_3}, we have
\begin{eqnarray}
\inf_{A\in \mathcal S_{\succeq}^n}{\|AX-B\|}_F^2 \hspace{-.3cm}&=& \hspace{-.3cm}\inf_{A_{11}\succeq 0,\,\operatorname{ker}(A_{11})\subseteq \operatorname{ker}(A_{21})}{\|A_{11}\Sigma_1-U_1^TBV_1\|}_F^2+ {\|A_{21}\Sigma_1-U_2^TBV_1\|}_F^2+{\|BV_2\|}_F^2 \nonumber \\
\hspace{-.3cm}&\geq &\hspace{-.3cm} \inf_{A_{11}\in \mathcal S_{\succeq }^r}{\|A_{11}\Sigma_1-U_1^TBV_1\|}_F^2+
\inf_{A_{21}\in \R^{n-r,r}}{\|A_{21}\Sigma_1-U_2^TBV_1\|}_F^2+{\|BV_2\|}_F^2
 \label{eq:proof_general_4}\\
\hspace{-.3cm}&=& \hspace{-.3cm}{\|\hat A_{11} \Sigma_1-U_1^TBV_1\|}_F^2 +{\|BV_2\|}_F^2, \label{eq:proof_general_5}
\end{eqnarray}
where the last equality follows by Lemma~\ref{lem:psd_fullrank} since (i) the first infimum
in the right hand side of~\eqref{eq:proof_general_4} is attained at a unique
$\hat A_{11} \in \mathcal S_{\succeq}^r$ (Lemma~\ref{lem:psd_fullrank}), and (ii) the value of the second infimum is zero, that is,
$\inf_{A_{21}\in \R^{n-r,r}}{\|A_{21}\Sigma_1-U_2^TBV_1\|}_F^2=0$, which is attained at $A_{21}=U_2^TBV_1\Sigma_1^{-1}$.
In order to show that equality in~\eqref{eq:proof_general_4} holds instead of inequality, we consider two cases.

{\underline{Case-1:}} When $\operatorname{ker}(\hat A_{11})\subseteq \operatorname{ker}(U_2^TBV_1\Sigma_1^{-1})$.
In this case, by taking  $A_{21}=U_2^TBV_1\Sigma_1^{-1}$ and $A_{11}=\hat A_{11}$ in~\eqref{eq:proof_general_1}, we define
\begin{equation}\label{eq:proof_general_6}
A_{opt}:=U\mat{cc}\hat A_{11} & (U_2^TBV_1\Sigma_1^{-1})^T\\ U_2^TBV_1\Sigma_1^{-1} & K\rix U^T,
\end{equation}
for some $K \in \R^{n-r,n-r}$ such that $K-(U_2^TBV_1\Sigma_1^{-1}){\hat A_{11}}^{\dag}(U_2^TBV_1\Sigma_1^{-1})^T \succeq 0$,
which upon simplification yields~\eqref{eq:psd_general_5}.
By Lemma~\ref{lem:psd_character}, we have  $A_{opt} \in \mathcal S_{\succeq}^n$ and,
in view of~\eqref{eq:proof_general_2}, it satisfies
\begin{equation}\label{eq:proof_general_7}
{\|A_{opt}X-B\|}_F^2={\|\hat A_{11}\Sigma_1-U_1^TBV_1\|}_F^2 + {\|BV_2\|}_F^2.
\end{equation}
This implies equality in~\eqref{eq:proof_general_4} in the case when
$\operatorname{ker}(\hat A_{11})\subseteq \operatorname{ker}(U_2^TBV_1\Sigma_1^{-1})$. This completes
the proof of part $\it{1}$.

{\underline{Case-2:}}
If suppose $\operatorname{ker}(\hat A_{11}) \nsubseteq \operatorname{ker}(U_2^TBV_1\Sigma_1^{-1})$, then let
\begin{equation}\label{eq:choose_eps}
0 < \epsilon < \left\{\begin{array}{ll}
\min\{1,{\|\hat A_{11}\Sigma_1-U_1^TBV_1\|}_F^2\} & \mbox{ if }\,
{\|\hat A_{11}\Sigma_1-U_1^TBV_1\|}_F \neq 0, \\
1 & \mbox{ otherwise, } \end{array}\right.
\end{equation}
and let $\hat A_{11}^\epsilon$ be as defined
in~\eqref{eq:psd_gen_6a}. Then we have
\begin{eqnarray*}
{\|\hat A_{11}^\epsilon \Sigma_1-U_1^TBV_1\|}_F^2 &=& {\|\hat A_{11}\Sigma_1-U_1^TBV_1 +
\hat U_2 \Upsilon \hat U_2^T \Sigma_1\|}_F^2 \nonumber \\
&\leq& {\|\hat A_{11} \Sigma_1-U_1^TBV_1\|}_F^2 + {\|\Upsilon\|}_F^2{\|\Sigma_1\|}_F^2+
2{\|\hat A_{11} \Sigma_1-U_1^TBV_1\|}_F{\|\Upsilon\|}_F{\|\Sigma_1\|}_F \nonumber\\
&<& {\|\hat A_{11} \Sigma_1-U_1^TBV_1\|}_F^2 +\epsilon,
\end{eqnarray*}
where the last inequality follows by using the fact that
$\epsilon$ satisfies~\eqref{eq:choose_eps},
%
and using the fact that
${\|\Upsilon\|}_F=\frac{\epsilon}{4{\|\Sigma_1\|}_F{\|\hat A_{11} \Sigma_1-U_1^TBV_1\|}_F}$
 when ${\|\hat A_{11} \Sigma_1-U_1^TBV_1\|}_F \neq 0$ and
 ${\|\Upsilon\|}_F=\frac{\epsilon}{4{\|\Sigma_1\|}_F}$, otherwise.
Thus we have
\begin{equation}\label{eq:proof_general_8}
{\|\hat A_{11}^\epsilon \Sigma_1-U_1^TBV_1\|}_F^2 < {\|\hat A_{11} \Sigma_1-U_1^TBV_1\|}_F^2 +\epsilon.
\end{equation}
Note that $\hat A_{11}^\epsilon$ is nonsingular. Thus again using $\hat A_{11}^\epsilon$ and by setting
$A_{21}=U_2^TBV_1\Sigma_1^{-1}$ in~\eqref{eq:proof_general_1}, we define
\begin{equation}\label{eq:for_rem}
A_{\epsilon}:=U\mat{cc}\hat A_{11}^\epsilon & (U_2^TBV_1\Sigma_1^{-1})^T\\ U_2^TBV_1\Sigma_1^{-1} & K_\epsilon\rix U^T,
\end{equation}
for some $K_\epsilon$ such that
$K_\epsilon-(U_2^TBV_1\Sigma_1^{-1})
{(\hat A_{11}^\epsilon)}^{-1}(U_2^TBV_1\Sigma_1^{-1})^T \succeq 0$. This upon
simplification yields~\eqref{eq:psd_general_6}, and, by construction, $A_{\epsilon} \in \mathcal S_{\succeq}^n$.
Thus in view of~\eqref{eq:proof_general_2} and~\eqref{eq:proof_general_8}, we have
\begin{eqnarray}\label{eq:proof_general_9}
{\|A_\epsilon X-B\|}_F^2 &=& {\|\hat A_{11}^\epsilon \Sigma_1-U_1^TBV_1\|}_F^2+{\|BV_2\|}_F^2 \nonumber \\
&<& {\|\hat A_{11} \Sigma_1-U_1^TBV_1\|}_F^2+{\|BV_2\|}_F^2 +\epsilon.
\end{eqnarray}
Thus as $\epsilon$ tends to zero, from~\eqref{eq:proof_general_5} and~\eqref{eq:proof_general_9}, we get the
equality in~\eqref{eq:proof_general_4}. Hence
\begin{equation*}
\inf_{A\in \mathcal S_{\succeq}^n}{\|AX-B\|}_F^2={\|\hat A_{11} \Sigma_1-U_1^TBV_1\|}_F^2+{\|BV_2\|}_F^2.
\end{equation*}
This infimum is attained when $\epsilon=0$, but then  by using Lemma~\ref{lem:psd_character} we have
$A_\epsilon \notin \mathcal S_{\succeq}^n$ because $\operatorname{ker}(\hat A_{11})
 \nsubseteq \operatorname{ker}(U_2^TBV_1\Sigma_1^{-1})$. Therefore the fact that
 $\inf_{A_{21}\in \R^{n-r,r}}{\|A_{21}\Sigma_1-U_2^TBV_1\|}_F^2=0$ and the uniqueness of $\hat A_{11}$
imply that the infimum is not attained. This completes the proof of $\it{2.}$
\eproof

\begin{remark}{\rm
As mentioned in Section~\ref{sec:contri}, part of Theorem~\ref{thm:general_analytic} has partially appeared in~\cite[Theorem 2.1]{Lia99}
under a differently stated necessary and sufficient condition, namely,
\[
\text{rank}(\hat A_{11})=\text{rank}([\hat A_{11} \, | \,   \Sigma_1^{-1}V_1^TB^TU_2]).
\]
However, our proof is different and constructive.
}
\end{remark}

\begin{remark}{\rm
In view of Theorem~\ref{thm:general_analytic}, a necessary and sufficient condition for the infimum
to be attained in~\eqref{psdproc} is that
$\operatorname{ker}(\hat A_{11}) \subseteq \operatorname{ker}(U_2^TBV_1\Sigma_1^{-1})$. This condition coincides with the result~\cite[Theorem 2.2]{Woo96}.
Therefore apart from the characterization of the solutions for the PSDP problem,
 Theorem~\ref{thm:general_analytic} also gives a completely different and relatively
 simpler proof for $\operatorname{ker}(\hat A_{11}) \subseteq \operatorname{ker}(U_2^TBV_1\Sigma_1^{-1})$
 to be a necessary and sufficient condition for the attainment of the infimum in~\eqref{psdproc}
 than the one provided in~\cite[Theorem 2.2]{Woo96}.

Note that when $\rank(X) = n$, $U_2$ is an $n$-by-$0$ empty matrix hence $\operatorname{ker}(U_2^T B V_1 \Sigma_1^{-1})$ is the full space so that the condition
$\operatorname{ker}(\hat A_{11}) \subseteq \operatorname{ker}(U_2^TBV_1\Sigma_1^{-1})$
for the infimum to be attained is always met.
Note also that when $\hat A_{11} \succ 0$ we have $\operatorname{ker}(\hat A_{11}) = \{0\}$ hence the condition $\operatorname{ker}(\hat A_{11}) \subseteq \operatorname{ker}(U_2^TBV_1\Sigma_1^{-1})$ is always met.
}
\end{remark}

\noindent

Using Lemma~\ref{lem:min_rank_F_2} in~\eqref{eq:proof_general_6} and in~\eqref{eq:for_rem} for matrices $A_{opt}$ and $A_\epsilon$, 
we can characterize the solutions of~\eqref{psdproc} with extremal properties of minimal
rank, minimal Frobenius norm or minimal spectral norm.

\begin{corollary}\label{cor:corollary1}
Let $X,B \in \R^{n,m}$, and let $r = \operatorname{rank}(X)$. Let also
$U_1,U_2,V_1,V_2$, $\Sigma_1$ and $\hat A_{11}$ be as in Theorem~\ref{thm:general_analytic}, and $Z := U_2^TBV_1\Sigma_1^{-1}$.
\begin{enumerate}
\item If $\operatorname{ker}(\hat A_{11}) \subseteq \operatorname{ker}(U_2^TBV_1\Sigma_1^{-1})$,
then $A_{opt}$ in~\eqref{eq:psd_general_5} with $K=Z\hat A_{11}^{\dag}Z^T$ is the unique solution of the problem~\eqref{eq:psd_general_1}
with minimal rank, minimal Frobenius norm and minimal spectral norm,
that is,
\[
\argmin_{A=\argmin_{A\in S^n_{\succeq}}{\|AX-B\|}_F}{\|A\|}_{F,2}~=~A_{opt}~=~
\argmin_{A=\argmin_{A\in S^n_{\succeq}}{\|AX-B\|}_F} {\rm rank}(A).
\]
\item Otherwise, for sufficiently small $\epsilon >0$, the matrix $A_{\epsilon}$ in~\eqref{eq:psd_general_6}
with $K_\epsilon=Z {(\hat A_{11}^\epsilon)}^{-1}Z^T$
is the unique matrix in $\mathcal S_\succeq^n$ with minimal rank, minimal Frobenius norm and minimal spectral norm, such that
\[
{\|A_\epsilon X-B\|}_F^2 < \inf_{A\in \mathcal S_{\succeq}^n}{\|AX-B\|}_F^2 + \epsilon.
\]
\end{enumerate}
\end{corollary}

\noindent
In the following theorem, we show that if $U_1^T(BX^T+XB^T)U_1 \preceq 0$ (with $U_1$ defined as in Theorem~\ref{thm:general_analytic})
then computing the exact value of the infimum in~\eqref{psdproc} does not require the solution $\hat A_{11}$ of the subproblem~\eqref{eq:subprob_A11} as in Theorem~\ref{thm:general_analytic}.
This complements the result in~\cite[Theorem 2.5]{Woo96} where the zero matrix is shown to be the unique solution of~\eqref{psdproc} when $\text{rank}(X)=n$ and $(BX^T+XB^T) \preceq 0$.

\begin{theorem}\label{thm:particular}
Let $X,B \in \R^{n,m}$, and let $r = \operatorname{rank}(X) < n$. Let also $U_1,U_2,V_1,V_2$ and $\Sigma_1$ be as defined in Theorem~\ref{thm:general_analytic}. If
$U_1^T(BX^T+XB^T)U_1 \preceq 0$, then
\begin{equation}\label{eq:psd_general_2}
\inf_{A\in \mathcal S_{\succeq}^n}{\|AX-B\|}_F^2={\|U_1^TBV_1\|}_F^2+{\|BV_2\|}_F^2,
\end{equation}
and it is not attained for any $A \in \mathcal S_{\succeq}^n$.
In this case, let
$\epsilon >0$ be sufficiently small and let $A_{11}^\epsilon \in \R^{r,r}$ be a diagonal matrix
with diagonal entries each equal to $\frac{\epsilon}{\alpha}$, where $\alpha=4\sqrt{n}{\|\Sigma_1\|}_F{\|U_1^TBV_1\|}_F$.
Define
\begin{equation}\label{eq:psd_general_3}
A_{\epsilon}:=U_1  A_{11}^\epsilon U_1^T + U_2(U_2^TBV_1\Sigma_1^{-1})U_1^T+ U_1 (U_2^TBV_1\Sigma_1^{-1})^T U_2^T
+ U_2 K_\epsilon U_2^T,
\end{equation}
where $K_\epsilon \in \R^{n-r,n-r}$ is such that
$K_\epsilon-(U_2^TBV_1\Sigma_1^{-1}) {(A_{11}^\epsilon)}^{-1}(U_2^TBV_1\Sigma_1^{-1})^T \succeq 0$.
Then $A_\epsilon \in \mathcal S_{\succeq}^n$ and
\begin{equation}\label{eq:psd_general_4}
{\|A_\epsilon X-B\|}_F^2 < \inf_{A\in \mathcal S_{\succeq}^n}{\|AX-B\|}_F^2  + \epsilon.
\end{equation}
\end{theorem}
\proof We only give a proof of~\eqref{eq:psd_general_2} and skip the proof
of~\eqref{eq:psd_general_4} as it is similar to the
proof of the point $\it{2}$ in Theorem~\ref{thm:general_analytic}.
In Theorem~\ref{thm:general_analytic} we proved that
\begin{equation}\label{eq:proof_general_10}
\inf_{A\in \mathcal S_{\succeq}^n}{\|AX-B\|}_F^2=
\min_{A_{11}\in \mathcal S_{\succeq}^r}{\| A_{11} \Sigma_1-U_1^TBV_1\|}_F^2+{\|BV_2\|}_F^2.
\end{equation}
Observe that for any $A_{11}\in \mathcal S_{\succeq}^r$, we have
\begin{eqnarray}\label{eq:proof_general_11}
{\| A_{11} \Sigma_1-U_1^TBV_1\|}_F^2 & =& \operatorname{trace}
\left((A_{11} \Sigma_1-U_1^TBV_1)^T(A_{11} \Sigma_1-U_1^TBV_1)\right) \nonumber \\
&=& {\|A_{11}\Sigma_1\|}_F^2 + {\|U_1^TBV_1\|}_F^2 -\operatorname{trace}\left(U_1^TBV_1\Sigma_1 A_{11}
+\Sigma_1V_1^TB^TU_1A_{11}\right) \nonumber \\
&=& {\|A_{11}\Sigma_1\|}_F^2 + {\|U_1^TBV_1\|}_F^2 - \operatorname{trace}\left(U_1^T(BX^T+XB^T)U_1A_{11}\right) \nonumber \\
&\geq& {\|U_1^TBV_1\|}_F^2,
\end{eqnarray}
where the third equality follows
by using the compact SVD of $X$, that is, $X=U_1\Sigma_1V_1^T$. The last inequality in~\eqref{eq:proof_general_11} follows by
using Lemma~\ref{lem:product_psd} since $\operatorname{trace}\left(U_1^T(BX^T+XB^T)U_1A_{11}\right) \leq 0$ as
$-U_1^T(BX^T+BX^T)U_1 \succeq 0$ and $A_{11}\succeq 0$. Therefore from~\eqref{eq:proof_general_11}, we obtain
\begin{equation*}
\min_{A_{11}\in \mathcal S_{\succeq}^r}{\| A_{11} \Sigma_1-U_1^TBV_1\|}_F^2 \geq {\|U_1^TBV_1\|}_F^2,
\end{equation*}
and equality holds when $A_{11}=0$. Plugging this in~\eqref{eq:proof_general_10}, we obtain
\begin{equation}\label{eq:proof_general_12}
\inf_{A\in \mathcal S_{\succeq}^n}{\|AX-B\|}_F^2={\|U_1^TBV_1\|}_F^2+{\|BV_2\|}_F^2.
\end{equation}
Using the arguments similar to that of Case-2 in Theorem~\ref{thm:general_analytic}, it follows that the infimum in~\eqref{eq:proof_general_12}
is not attained and for a sufficiently small $\epsilon > 0$, $A_\epsilon$ in~\eqref{eq:psd_general_3}
satisfies~\eqref{eq:psd_general_4}. This completes the proof.
\eproof

\begin{remark}{\rm
Theorem~\ref{thm:general_analytic} reduces the original problem
$(\mathcal P)$ to a PSDP problem with a diagonal $r$-by-$r$ matrix $X$, with $r=\rank(X)$,  of the form
\begin{equation}\label{eq:temp1}
\min_{\widetilde A \succeq 0}{\|\widetilde A \Sigma -\widetilde B\|}_F^2,
\end{equation}
where $\Sigma \in \mathcal{S}^{r}_{\succ}$ is a diagonal matrix with positive diagonal entries, for which a unique solution
is guaranteed by Lemma~\ref{lem:psd_fullrank}.
In some special cases, when $ \widetilde B^T\Sigma + \Sigma^T \widetilde B \succeq 0$~\cite{Lia99},
or when $ -\widetilde B \Sigma^T-\Sigma \widetilde B^T \succ 0$~\cite{Woo96},
the optimal solution in~\eqref{eq:temp1} can be explicitly given.
However, in general finding an analytic solution to the subproblem is a challenging task and still
an open problem. We will discuss in Section~\ref{fgm} the use of an optimal first-order method to solve this problem.
}
\end{remark}

\subsection{Computational cost of the semi-analytical approach} \label{companal}

In view of Theorem~\ref{thm:general_analytic}, we have that the semi-analytical approach completes in three steps.
The first step takes $O(\max(m,n)\min(m,n)^2)$ floating point operations
to compute the singular value decomposition of $X$ ~\cite{TreB97}.
The second step is to compute the solution $\hat A_{11}$ of the subproblem~\eqref{eq:subprob_A11},
and cost of this depends on the method used to solve it (see Section~\ref{fgm}).
The third step is to form $A_{opt}$. This involves matrix-matrix multiplication and costs $O(n^2 r)$.

\subsection{ Analytical solution for $\rank(X) = 1$ }\label{sec:rank_1}
We note that the case $\rank(X) = 1$ differs from the case $\rank(X) > 1$ because in that case the subproblem~\eqref{eq:subprob_A11} has a closed-form solution:
\[
\hat A_{11} = \max \left( 0 , U_1^TBV_1 \Sigma_1^{-1} \right) \; \in \;  \mathbb{R}.
\]
Therefore, we can provide a complete analytical characterization of the set of optimal solution of~\eqref{psdproc} in this particular case. Although this follows directly from Theorem~\ref{thm:general_analytic} and the observation above,
we state the result here in the case  $\operatorname{rank}(X)=1$  for the sake of future reference.

The rank-one case is also important in solving the particular case when $X$ and $B$ are vectors, that is, when one is looking for $A \succeq 0$ such that $\| A x - b \|_2$ is minimized where $x$ and $b$ are vectors.

\begin{theorem}\label{thm:analytic_rank1}
Let $X,B\in \R^{n,m}$ be such that $\operatorname{rank}(X)=1$.
Let $X=U\Sigma V^T$ be a singular value decomposition
of $X$, where $U=[u~\,U_1] \in \R^{n,n}$ with $u \in \R^n$, $V=[v~\,V_1]\in \R^{m,m}$ with $v \in \R^m$, and
$\Sigma=\mat{cc}\sigma & 0 \\0& 0\rix \in \R^{n,m}$ with $\sigma > 0$.
Then the following hold.
\begin{enumerate}
\item If $u^TBv > 0$, then
\[
\inf_{A \in \mathcal S_{\succeq}^n} {\|AX-B\|}_F = {\|BV_1\|}_F,
\]
and $A_{opt}$ attains the infimum 
if and only if
%
\begin{equation}\label{eq:rank_1_result_1}
A_{opt}=\sigma^{-1}\left((u^TBv)uu^T + U_1U_1^TBvu^T + uv^TB^TU_1U_1^T \right) + U_1 K U_1^T,
\end{equation}
for some matrix $K$ such that $K-\frac{1}{\sigma u^TBv}U_1^T (Bv)(Bv)^T U_1 \succeq 0$.
In particular, $ A_{opt}$ can be chosen to be of rank one by choosing $K=\frac{1}{\sigma u^TBv}U_1^T (Bv)(Bv)^T U_1$.
\item If $u^TBv \leq 0$, then
\begin{equation}\label{eq:rank_1_result_2}
\inf_{A \in \mathcal S_{\succeq}^n} {\|AX-B\|}_F^2={\|u^TBv\|}_F^2 + {\|BV_1\|}_F^2.
\end{equation}
Further, if $U_1^TBv=0$, then the infimum in~\eqref{eq:rank_1_result_2} is attained by a
matrix $A_{opt}$ of the form~\eqref{eq:rank_1_result_1}.
If $U_1^TBv \neq 0$, then the infimum in~\eqref{eq:rank_1_result_2} is not attained. In the later case
for any arbitrary small $\epsilon > 0$, choose $n_0 \in \mathbb N$  such that
$\frac{\sigma^2}{n_0^2}-2\sigma\frac{u^TBv}{n_0} < \epsilon$
and define
\[
 A_{n_0}=\sigma^{-1}\left(\frac{1}{n_0} uu^T + U_1U_1^TBvu^T + uv^TB^TU_1U_1^T\right)  + U_1  K_{n_0} U_1^T,
\]
for some $K_{n_0}$ with $K_{n_0}-\frac{n_0}{\sigma^2}U_1^T (Bv)(Bv)^T U_1 \succeq 0$. Then  $A_{n_0} \succeq 0$ and
\[
{\|A_{n_0}X-B\|}_F^2 < \inf_{A\succeq 0}{\|AX-B\|}_F^2 + \epsilon.
\]
Moreover, $A_{n_0}$ can be chosen to be of rank one by choosing
$K_{n_0}=\frac{n_{0}}{\sigma^2 }U_1^T (Bv)(Bv)^T U_1$.
\end{enumerate}
\end{theorem}

\section{An algorithmic solution to the PSDP problem}

In this section, we first describe the fast gradient method to solve~\eqref{psdproc} (Section~\ref{fgm}).
We then propose a new very efficient initialization strategy for~\eqref{psdproc} when $X$ is diagonal and ill-conditioned (Section~\ref{init}).
Finally, we explain the advantages of combining the semi-analytical approach, FGM and our new initialization strategy.
In particular, this allows us to guarantee linear convergence for solving~\eqref{psdproc} (Section~\ref{ultimateapproach}).

\subsection{Fast Gradient Method} \label{fgm}
In order to be able to solve large-scale problem~\eqref{psdproc}, say for $n$ up to a 1000,
it makes sense to use first-order methods.
In this section, we describe the fast gradient method (FGM) applied to~\eqref{psdproc};
see Algorithm~\ref{fastgrad}. We choose FGM because it is simple to implement and it is an optimal first-order method for smooth convex optimization, that is, no first-order method can have a faster convergence rate.
In fact, FGM is guaranteed to decrease the objective function value at sublinear rate $O(1/t^2)$ where $t$ is the iteration number. Moreover, in the strongly convex case, when $\kappa = \frac{\sigma_1(X)}{\sigma_n(X)} > 0$,
the decrease is guaranteed to be at linear rate $O( (1-1/\kappa)^t )$.
This is much faster than the standard gradient descent method, with respective rate of $O(1/t)$ and $O( (1-1/\kappa^{2})^t )$.
We refer the reader to~\cite{nes83, nes04} for more details on FGM.

\renewcommand{\thealgorithm}{FGM}
 \algsetup{indent=2em}
\begin{algorithm}[ht!]
\caption{Fast Gradient Method for~\eqref{psdproc} \cite[p.90]{nes04}} \label{fastgrad}
\begin{algorithmic}[1]
\REQUIRE An initial guess $A \in \mathcal{S}^n_{\succeq}$,
number of iterations $T$ (other stopping criteria can be used).
\ENSURE An approximate solution $A \approx \argmin_{\tilde{A} \succeq 0} \|\tilde{A} X - B\|_F$.  \medskip
\STATE $L = \sigma_{1}^2(X)$, $q = \frac{\sigma_{n}^2(X)}{L}$.
\STATE $\alpha_1 \in (0,1)$; $Y = A$.
\FOR{$k = 1 :$ $T$}
\STATE $\hat A = A$. \hspace{2.4cm} \emph{\% Keep the previous iterate in memory.}
\STATE $G_Y = Y X X^T - B X^T$. \emph{\% Compute the gradient, $XX^T$ and $BX^T$ can be pre-computed.}
\STATE $A = \mathcal{P}_{\succeq 0}  \left(  Y - \frac{1}{L} G_Y \right)$. \hspace{0.05cm} \emph{\% Projected gradient step.}
\STATE $\alpha_{k+1} = \frac{1}{2} \left( q - \alpha_{k}^2  +\sqrt{ (q - \alpha_{k}^2 )^2 +4 \alpha_{k}^2 } \right)$, $\beta_k =  \frac{\alpha_{k} (1-\alpha_{k})}{\alpha_{k}^2 + \alpha_{k+1}}$.
\STATE $Y = A + \beta_k \left(A - \hat A\right)$.  \hspace{0.02cm} \emph{\% Linear combination of the current and previous iterates.}
\ENDFOR
\end{algorithmic}
\end{algorithm}

The computational cost of FGM is $O(n^3 + n^2 m)$ operations per iteration.  The most expensive steps are
\begin{itemize}
\item the computation of the singular values of $X$ (step 1) requiring $O\left(\min(m,n)^2 \max(m,n)\right)$
operations~\cite{TreB97}.
\item the computation of the gradient in $O(mn^2 + n^3)$ operations (step 5).
Note that the $n$-by-$n$ matrices $XX^T$ and $BX^T$ should be computed only once in which case the remaining computational cost for the gradient computation per iteration is $O(n^3)$.
\item The projection step in $O(n^3)$ operations (step 6), as it requires the eigenvalue decomposition of a symmetric $n$-by-$n$ matrix;
see~\eqref{projpsd}.
\end{itemize}
Denoting $T$ the total number of iterations (typically, $T \geq 100$), the total computational cost of FGM
is $O(T n^3 + mn^2)$. In most cases, $Tn \geq m$ hence the computational cost of FGM will be $O(T n^3)$.

\subsection{Initialization} \label{init}
In this section, we present three initialization strategies.

\subsubsection{Projection of the optimal unconstrained solution}  \label{uncstr}
In~\cite{AE97}, the authors propose to use as an initialization the projection of the optimal solution of the unconstrained problem, that is,
\[
\mathcal{P}_{\succeq} \big( \argmin_{A \in \mathbb{R}^{n \times n}} \|AX - B\|_F \big) .
\]
This initialization can sometimes perform well. However, it comes with no guarantee and provides very bad initialization in several situations,
in particular for ill-conditioned problems; see Section~\ref{prelimNE} below for some examples.

\subsubsection{Diagonal matrix}  \label{diag}
It is rather straightforward to compute the optimal solution of~\eqref{psdproc} assuming that the matrix $A$ is diagonal.
In fact, the problem reduces to $n$ independent least squares problem in one variables with a nonnegativity constraint. We have
\[
 \left\|\diag(a) X - B\right\|_F^2 =  \sum_{i=1}^n { \left\|a_i X(i,:) - B(i,:)\right\|}_F^2,
\]
where $\diag(a)$ is a diagonal matrix whose diagonal elements are given by
$a=[a_1,a_2,\ldots,a_n]^T \in \mathbb{R}^n$.
The optimal solution for each subproblem is given by
\[
a_i^* =
\argmin_{a_i \geq 0} \left\|a_i X(i,:) - B(i,:)\right\|_F^2  =
\max \left( 0, \frac{B(i,:)X(i,:)^T }{{\| X(i,:) \|}_2^2} \right), \quad i=1,2,\dots n.
\]

\subsubsection{Recursive decomposition for ill-conditioned and diagonal $X$}  \label{recinit}

As we have explained previously, the convergence of first-order algorithms for~\eqref{psdproc} will depend on the conditioning of $X$.
Using the semi-analytical approach from Section~\ref{semianalsol}, \eqref{psdproc} can be reduced to a problem where $X$ is diagonal with positive diagonal elements.
If $X$ is well-conditioned, then FGM will converge fast and the initialization strategy does not play a crucial role.
However, when $X$ is ill-conditioned, FGM will be more sensitive to initialization
as it converges slower.

Let us generalize the idea from the previous section by assuming that $A$ is block diagonal instead of diagonal. For simplicity, let us assume that the diagonal matrix $X$  is partitioned into to blocks $X_1$ and $X_2$
(this generalizes easily to more than two blocks):
\[
\min_{A_{1} \succeq 0, A_{2} \succeq 0}
{\left\|
\left( \begin{array}{cc} A_{1} & 0 \\ 0 & A_{2} \end{array} \right)
\left( \begin{array}{cc} X_{1} & 0 \\ 0 & X_{2} \end{array} \right)
- \left( \begin{array}{cc} B_{1} & B_{12} \\ B_{21} & B_{2} \end{array} \right)   \right\|}_F .
\]
This problem can be decoupled into two independent subproblems: for $i=1,2$,
\[
\min_{A_{i} \succeq 0} {\|A_{i} X_{i} - B_{i}\|}_F.
\]
Let us denote $\kappa(X_{i})$ the condition number of $X_{i}$.
If $X_1$ and $X_2$ are well conditioned, that is, $\max_i \kappa(X_{i})$ is small,
good approximate solutions to these subproblems can be obtained much faster than for the ill-conditioned $X$.
%
Moreover, since $X$ is diagonal, partitioning $X$ into two blocks in order to minimize $\max_i \kappa(X_i)$ can be done as follows:
\begin{enumerate}
\item Sort the diagonal entries of $X$ such that $x_1 \leq x_2 \leq \dots \leq x_n$
(in our reduction, the diagonal entries of $X$ are already sorted in nonincreasing order,
since we use the standard SVD). This can be done in $O(n \log(n))$ operations.
\item Pick the partition $[1, 2, \dots, k] \cup [k+1, k+2, \dots, n]$ such that $\max \left( \frac{x_k}{x_1}, \frac{x_n}{x_{k+1}} \right)$ is minimized.  This can be done in $O(n)$ operations.
\end{enumerate}
Finally, it is straightforward to use this idea recursively as follows:
As long as a block $X_i$ is not well-conditioned, that is, $\kappa(X_i) > \kappa_M$ for some parameter $\kappa_M$ (we use $\kappa_M =100$),
the block $X_i$ is partitioned into two blocks as explained above.
Once $X$ has been partitioned into well-conditioned subblocks, we combine the diagonal initialization along with 100 iterations of FGM to approximately solve the (well-conditioned) subproblems.
As we will see, not only this initialization will provide a solution with low initial error,
it will also allow the FGM to converge faster to the optimal solution.

Note that this initialization is only applicable when $X$ is diagonal hence can be used only
in combination with the semi-analytical reduction described in Section~\ref{semianalsol}.

\subsubsection{Preliminary numerical experiments} \label{prelimNE}
\begin{figure*}[h!]
\begin{center}
\caption{Evolution of the average error $\|AX - B\|_F$ over 100 randomly generated matrices $B$ for the different initializations using FGM: Gaussian (left) and Uniform (right).
\label{figinit}}
\begin{tabular}{cc}
\includegraphics[width=0.45\textwidth]{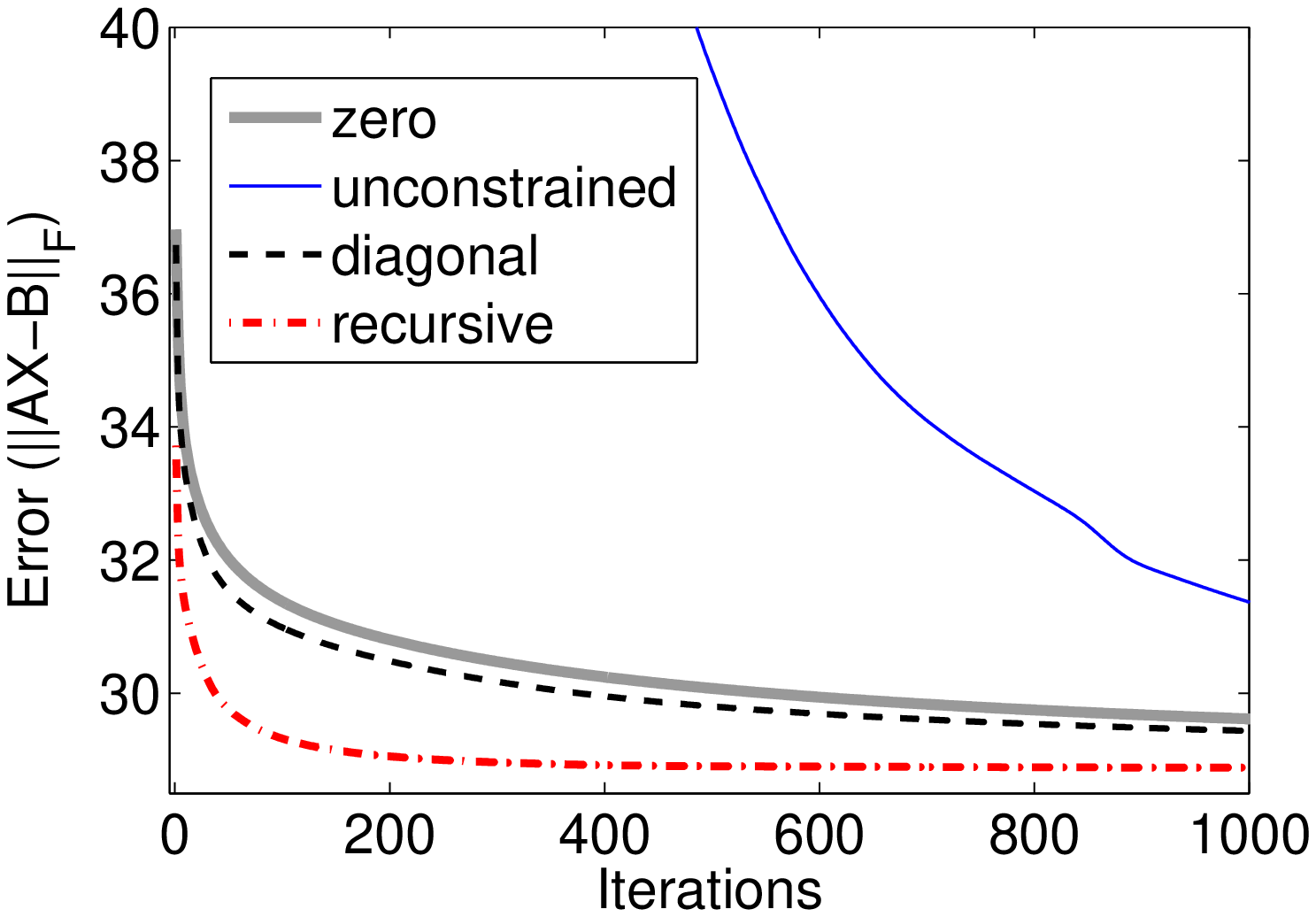}  & \includegraphics[width=0.45\textwidth]{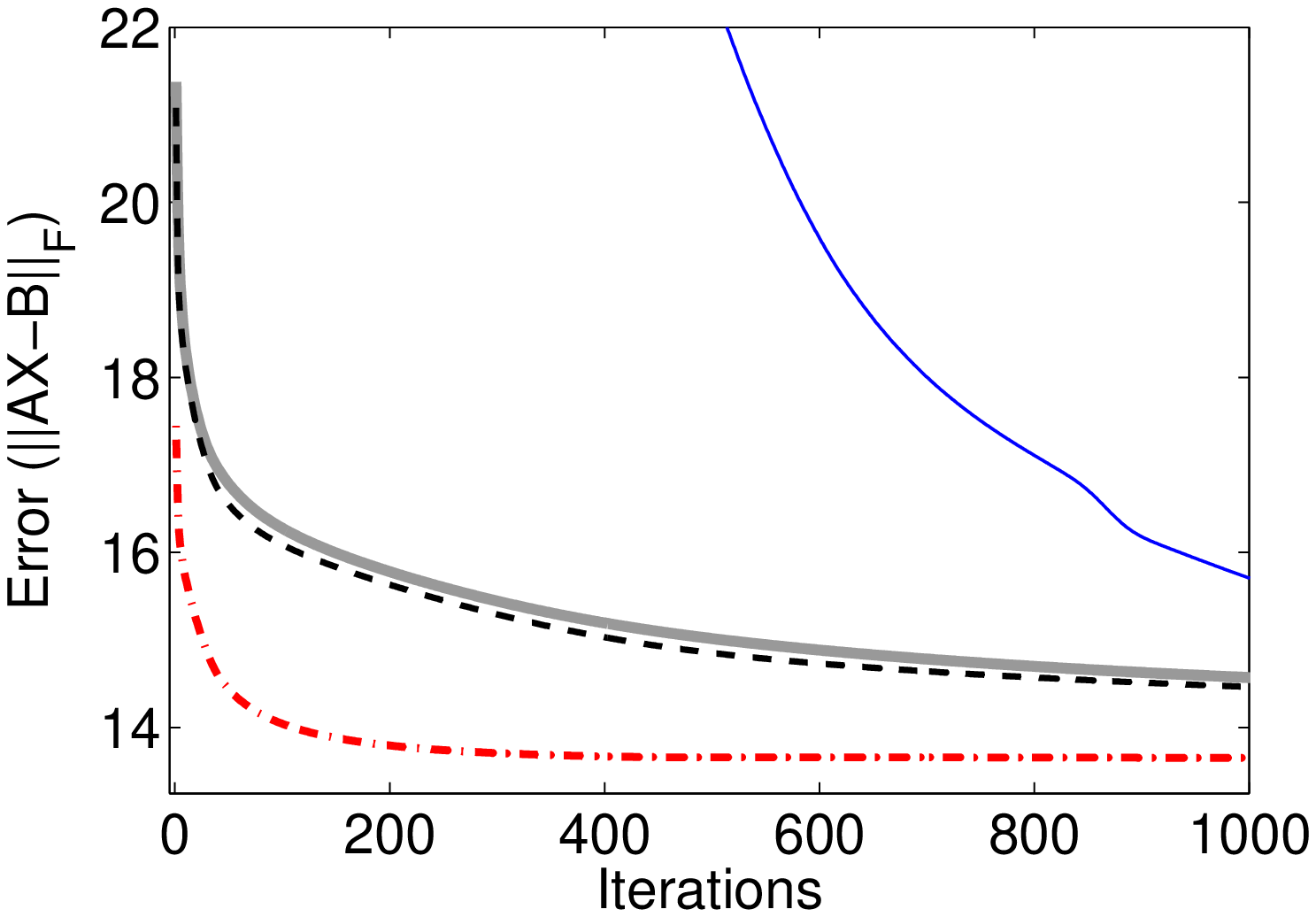}
\end{tabular}
\end{center}
\end{figure*}

Let us compare the different initialization strategies in ill-conditioned cases (well-conditioned cases are not so interesting since most algorithms will converge fast, being less sensitive to initialization),
and let us consider
\[
X = \diag(1,2,\dots,10, 20, \dots 100, 200 \dots, 1000, 2000, \dots, 10000), \text { with } \kappa(X) = 10^4.
\]
Recall that the optimal solution of~\eqref{psdproc} will be unique in this case since $X$ has rank $n=37$ (Lemma~\ref{lem:psd_fullrank}).
We generate $B$ in two different ways:
\begin{itemize}
\item Gaussian. Each entry is randomly generated following a normal distribution with mean 0 and standard deviation 1 (\texttt{randn(n)} in Matlab).
\item Uniform. Each entry is randomly generated following a uniform distribution in the interval [0,1] (\texttt{rand(n)} in Matlab).
\end{itemize}
In each case, we generate 100 such matrices.  We will compare four initializations:
\begin{itemize}
\item Zero. This is the trivial initialization $A = 0$.
\item Unconstrained. This is the projection onto $\mathcal{S}^n_{\succeq}$ of the optimal solution of the unconstrained version of~\eqref{psdproc}; see Section~\ref{uncstr}.
\item Diagonal. This is the optimal diagonal solution of~\eqref{psdproc}; see Section~\ref{diag}.
\item Recursive. This is the recursive decomposition approach described in Section~\ref{recinit}.
\end{itemize}
We will only compare the initializations combined with Algorithm~\ref{fastgrad} because, as we will see in Section~\ref{numexp}, it consistently performs well.
Table~\ref{figinit} gives the average initial error for the four initialization strategies, along with the standard deviation, for the 100 randomly generated matrices $B$ of the two types.
The best result is highlighted in bold.
Figure~\ref{tableinit} displays the evolution of the error (average over 100 runs) for the different initializations using FGM.

\begin{center}
 \begin{table}[h!]
 \begin{center}
\caption{Average initial error ${\|AX - B\|}_F$ and standard deviation (in brackets)
obtained by the different initialization approaches.
\label{tableinit}}
 \begin{tabular}{|c|c|c|c|c|}
 \hline   &  Zero & Unconstrained & Diagonal & Recursive   \\
 \hline
Gaussian &  36.97 (0.81)  &  8764.30 (1673.09)   &  36.73 (0.81)   & \textbf{33.72} (0.78)   \\ \hline
Uniform &  21.37 (0.27)  &  5799.61 (726.11)   &   21.08 (0.26)   & \textbf{17.45} (0.29)   \\ \hline
\end{tabular}
 \end{center}
 \end{table}
 \end{center}

We observe the following
\begin{itemize}
\item The initialization based on the projection of the optimal solution of the unconstrained problem performs very badly.
In fact, the error is significantly larger than with the trivial zero initialization.
\item The diagonal initialization performs slightly better than the zero initialization.
\item Our recursive initialization performs best,
both in term of initial error and for enabling FGM to converge faster to the unique optimal solution of~\eqref{psdproc}.
\end{itemize}

\subsection{ Combination of the semi-analytical approach, the recursive initialization and FGM } \label{ultimateapproach}

Combining the semi-analytical approach, the recursive initialization and FGM allow us to obtain an efficient algorithm for~\eqref{psdproc}.
The semi-analytical approach reduces the problem~\eqref{psdproc} to a problem
(i) involving only (possibly smaller) square matrices (where $m=n=\rank(X)$) and
(ii) that is strongly convex (since the `new' $X$, denoted $\Sigma_{1}$ in our derivations,
is diagonal with positive diagonal elements).
This requires $O(\min(m,n)^2 \max(m,n) + n^2 r)$ operations; see Section~\ref{companal}.

The first advantage (i) of the semi-analytical approach is that it significantly reduces the computational cost of FGM when $\rank(X) \ll n$ (for example when $m \ll n$).
The second advantage (ii) guarantees FGM to decrease the objective function value at linear rate
$( 1 - 1/\kappa_r )$ where $\kappa_r = \frac{\sigma_1(X)}{\sigma_r(X)}$ and $r = \rank(X)$.
This is, to the best of our knowledge, the first time an algorithm is proposed for~\eqref{psdproc} with guaranteed linear convergence.
If $X$ is ill-conditioned, that is, $\kappa_r$ is large, the convergence could be slow. However, this is mitigated by our recursive decomposition strategy that solves well-conditioned subproblems to initialize the ill-conditioned one.

\section{Numerical Experiments} \label{numexp}

In this section, we compare the following algorithms:
\begin{itemize}
\item Gradient. The projected gradient method applied on~\eqref{psdproc} --this is FGM using $\beta_k = 0$ at each step (that is, $A=Y$ in Algorithm FGM). This method will serve as a baseline.
\item FGM. The fast projected gradient method applied on~\eqref{psdproc}; see Algorithm~\ref{fastgrad}.
\item ParTan. This is the method proposed in~\cite{AE97} and referred to as `Parallel tangents'.
This algorithm is rather similar to FGM but does not guarantee the optimal convergence rate.
It can be seen as a heuristic variant of FGM where $\beta_k$ is chosen as to minimize $\|YX-B\|_F$ without the PSD constraint,
where $Y = A + \beta_k \hat A$ with $A$ the current iterate and $\hat A$ the previous iterate; see Algorithm~\ref{fastgrad}.
If this step does not decrease the objective function, then $\beta_k = 0$ is chosen (that is, a standard gradient step is used). Note that the computation of the $\beta_k's$ makes ParTan computationally slightly more expensive than FGM and Gradient.
\item AN-FGM. This is the combination of the semi-analytical approach,
reducing the problem to the case where $X$ is diagonal with positive diagonal elements, and then using FGM on this reduced problem.
We use the recursive initialization described in Section~\ref{recinit}.
\end{itemize}
Note that we could combine the analytical approach and the recursive initialization with any other method.
We choose FGM because it guarantees linear convergence, although it performs similarly as ParTan (see the numerical experiments below).
For the first three algorithms, we use the diagonal initialization.

In these numerical experiments, we try as much as possible to cover all the different scenarios:
we test for $m=n$, $m < n = 2m$ and $n < m=2n$. In all cases, the matrix $B$ is generated in the same way: each entry is randomly generated following a normal distribution with mean 0 and standard deviation 1 (\texttt{randn(m,n)} in Matlab).
For the matrix $X$, we consider three cases

\begin{enumerate}
\item Well-conditioned. Each entry is randomly generated following a normal distribution with mean 0 and standard deviation 1 (\texttt{randn(m,n)} in Matlab).
\item Ill-conditioned. Let $(U,\Sigma,V)$ be the compact SVD of a matrix generated as in the well-conditioned case. Then we generate $X = U \Lambda V$ where $\Lambda$ is a diagonal matrix such that $\Lambda(i,i) = \alpha^{i-1}$ and $\alpha^{\min(m,n)-1} = 10^6 = \kappa(X)$.
\item Rank deficient. We perform the SVD $(U,\Sigma,V)$ of a matrix generated as in the well-conditioned case, set the $r = \min(m,n)/2$ smallest singular values
of $\Sigma$ to zero to obtain $\Sigma'$, and then compute $X = U \Sigma' V^T$ so that $\rank(X) = \min(m,n)/2$.
\end{enumerate}

The Matlab code is available from \url{https://sites.google.com/site/nicolasgillis/}.
All tests are preformed using Matlab on a laptop Intel CORE i5-3210M CPU @2.5GHz 2.5GHz 6Go RAM.

For each experiment, we generate 10 such matrices and display the average results. For each algorithm, we perform 1000 iterations.

\begin{figure}[h!]
\begin{center}
\caption{Evolution of the average relative error $100 \frac{\|AX-B\|_F}{\|B\|_F}$
for the different algorithms.
\label{resultsALL}}
\includegraphics[width=\textwidth]{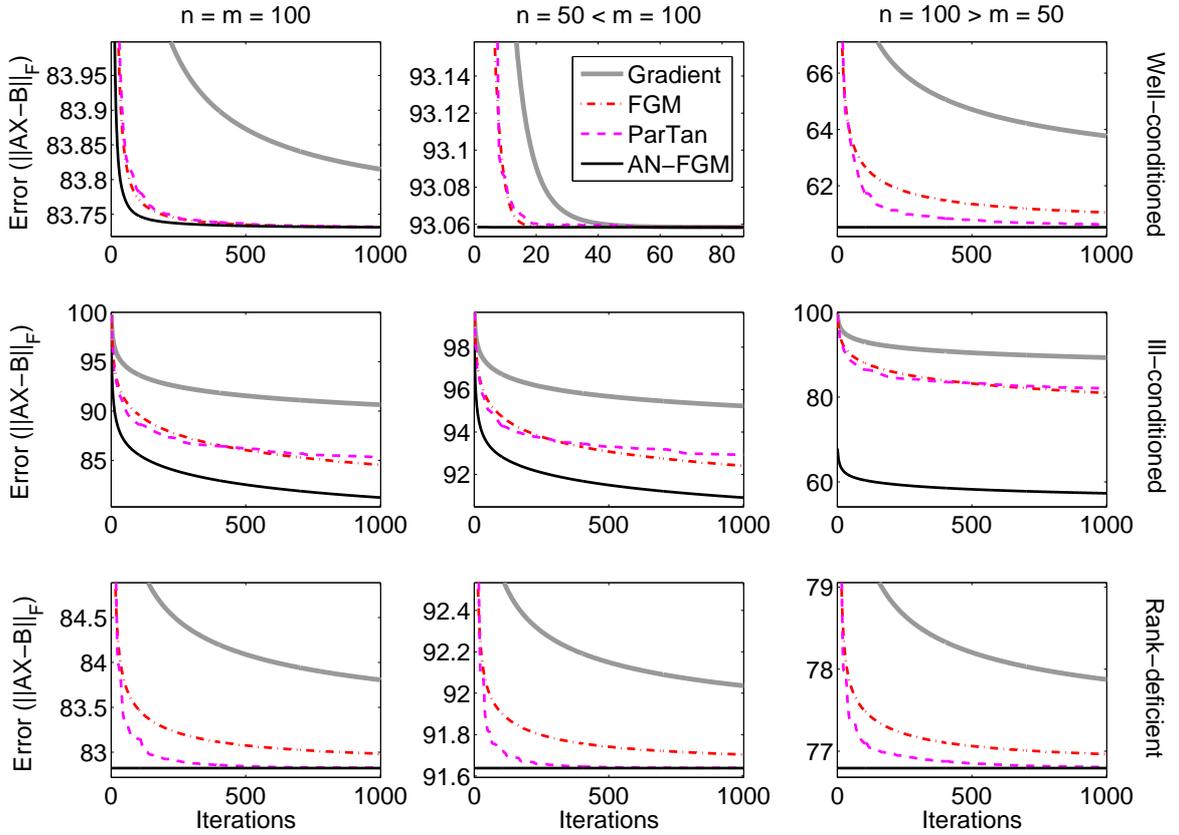}
\end{center}
\end{figure}

Figure~\ref{resultsALL} displays the evolution of the relative error in percent, that is,
\[
\text{relative error (\%)} = 100 \frac{\|AX-B\|_F}{\|B\|_F} ,
\]
for each algorithm in each of the nine cases described above. Table~\ref{tab:results} reports the computational time required by each algorithm to perform the 1000 iterations.
In bold we indicate the cases when AN-FGM  is significantly faster than the other approaches,
because of the dimension reduction of the problem.

\begin{center}
 \begin{table}[h!]
 \begin{center}
\caption{Computational times for the different algorithms to perform 1000 iterations.
\label{tab:results}}
 \begin{tabular}{|c|c||c|c|c|c|}
 \hline             &           &  Gradient &  FGM  &  ParTan  &  AN-FGM
 \\ \hline
 Well-conditioned & $m=n=100$  & 6.71   &  6.57 & 8.16 & 7.01 \\
                  & $m=2n=100$ & 1.48   &  1.49 & 2.13 & 1.52 \\
                  & $n=2m=100$ & 5.58   &  5.59 & 6.31 & \textbf{1.45} \\ \hline
 Ill-conditioned & $m=n=100$  & 5.87   &  5.67 & 7.49 & 5.64 \\
                  & $m=2n=100$ & 1.37   &  1.32 & 2.00 & 1.33 \\
                  & $n=2m=100$ & 5.49   &  5.20 & 6.75 & \textbf{1.28} \\ \hline
 Rank-deficient   & $m=n=100$  & 6.22   &  6.56 & 7.27 & \textbf{1.58} \\
                  & $m=2n=100$ & 1.67   &  1.42 & 2.22 & \textbf{0.34} \\
                  & $n=2m=100$ & 5.33   &  5.33 & 6.55 & \textbf{0.36} \\ \hline
\end{tabular}
 \end{center}
 \end{table}
 \end{center}

We observe the following:
\begin{itemize}
\item In terms of computational time, AN-FGM will be faster when $\rank(X) \ll n$
since, after the preprocessing performed by our semi-analytical approach, the number of operations per iteration of AN-FGM
is $O(r^3)$ where $r = \rank(X) \leq n$. This happens when $X$ is rank-deficient, and when $m = n/2$ --these are the bold results in Table~\ref{tab:results}. In all other cases, all algorithms have roughly the same computational cost, namely $O(Tn^3)$ where $T$ is the number of iterations (here $T=1000$); see the discussion in Section~\ref{fgm}.
\item In all cases, the gradient method performs the worse. This is not surprising since it only uses the gradient information of the current iterate as opposed to FGM and ParTan.
\item For well-conditioned $X$:

For $m=n=100$, FGM, ParTan and AN-FGM perform similarly. The reason is that the
semi-analytical approach cannot reduce the problem. The slight advantage of AN-FGM compared to FGM and ParTan comes from the fact the the recursive initialization already computes 100 iterations of FGM.

For $n = 50 < m = 100$, all algorithms converge very fast, within 100 iterations:
the reason is that the condition number of $X$ is much smaller than in the case $n = m = 100$
(the average condition number of a 50-by-100 Gaussian matrix is below 10, while for a 100-by-100 Gaussian matrix it is above 1000).

For $n = 100 > m = 50$,  AN-FGM performs best because it reduces the dimension of the problem from $100$ to $50$ (hence also reducing the computational cost per iteration; see Table~\ref{tab:results}).

\item For ill-conditioned $X$, AN-FGM outperforms the other approaches, especially for $n=2m = 100$ (for the same reasons as in the well-conditioned case). For $m=n=100$ and $n = 50 < m = 100$, the better performance of AN-FGM is explained by the recursive initialization; see Section~\ref{prelimNE}. FGM and ParTan perform similarly, with a slight advantage for FGM.
\item For rank-deficient $X$, it is easy to analyze: AN-FGM outperforms the other approaches because
 it reduces the problem size (from an $n$-by-$n$ variable problem to an $r$-by-$r$ where $r \ll n$) and the reduced problem is well-conditioned (since the non-zero singular values of $X$ come from a randomly generated matrix).
It is interesting and surprising to note that ParTan performs better than FGM in this case\footnote{ Although we observed that by tuning the parameter $\alpha_0$
(namely, using 0.9 instead of 0.1) makes FGM perform slightly better than ParTan in these rank-deficient cases. It could be an interesting direction for research to tune this parameter automatically. }.
However, we prefer not to use ParTan because it is a heuristic to combine several iterates and comes with no guarantee (authors only prove that there is at least one subsequence of the iterates converging to the optimal solution~\cite[Lemma~4.1]{AE97}).
\end{itemize}

\noindent \paragraph{Summary} 
In all cases, the gradient method performs significantly worse than the other first-order methods.
When $r = \rank(X) \ll n$ or $X$ is ill-conditioned, AN-FGM outperforms FGM and ParTan.
In fact, when $r = \rank(X) \ll n$, each iteration of AN-FGM is in $O(r^3)$ operations instead of $O(n^3)$,
and, when $X$ is ill-conditioned, AN-FGM takes advantage of an effective initialization strategy.
In the other cases (that is, $\rank(X) \approx n$ and $X$ is well-conditioned),
AN-FGM, FGM and ParTan perform similarly.

\begin{remark}{\rm \label{reminfnotat}
For all generated matrices in the rank-deficient cases and the cases $n = 2m$, we have $\rank(X) < n$ implying that the infimum of~\eqref{psdproc} is not necessarily attained; see Lemma~\ref{lem:psd_fullrank}.
In fact, we have observed that the infimum is never attained.
Although we do not have a rigorous explanation for this fact,
we believe that for random matrices $X$ and $B$ with $\rank(X) < n$,
 it is very unlikely that the condition
 $\operatorname{ker}(\hat A_{11}) \subseteq \operatorname{ker}(U_2^T B V_1 \Sigma_1^{-1})$
 for the infimum to be attained is met. Recall that the matrix $\hat A_{11}$ is the solution of the
 PSD Procrustes subproblem~\eqref{eq:subprob_A11}, and $U_2$, $V_1$ and $\Sigma_1^{-1}$ are factors in the SVD of $X$;
see Theorem~\ref{thm:general_analytic}. We believe these conditions are not likely to be met when $n$ is large because
\begin{itemize}
\item  The solution $\hat A_{11}$ of the subproblem~\eqref{eq:subprob_A11} is in general not positive definite.
Take for example the simple case $\Sigma_1 = I$ for which $\hat A_{11}$
is the projection of $(C+C^T)/2$ on the cone of PSD matrices where $C = U_1^T B V_1$. Since $U_1$ and $V_1$ come from the SVD of $X$ which is randomly generated, and $B$ is randomly generated, the entries of $C$ also follows a Gaussian-like distribution (for which the probability for an eigenvalue to be positive is 1/2).
Therefore,  it is not likely for $(C+C^T)/2$ to be positive definite hence its projection is in most cases rank deficient; see~\eqref{projpsd}.
\item The kernel condition is not likely to be satisfied: the probability for one subspace to contain another subspace generated randomly is zero (of course, in our case the subspaces are not independent so a rigorous probabilistic analysis is non trivial).
\end{itemize}
}
\end{remark}

\subsection{Comparison with a second-order method}

In this section, we compare AN-FGM with the interior point method SDPT3 (version 4.0)~\cite{toh1999sdpt3, tutuncu2003solving}, where we used CVX as a modeling system~\cite{cvx, gb08}. This is a second-order method hence it is computationally more expensive but guarantees quadratic convergence.
We perform in this section exactly the same experiment as in the previous section except that we use matrices of smaller size (for the sizes of the previous section, SDPT3 needs more than one minute to terminate).
In order to have a fair comparison, we first run SDPT3 on~\eqref{psdproc} and then run AN-FGM allowing the same computational time as for SDPT3.
Table~\ref{tab:resultsIPM} gives the relative error in percent for the different types of matrices for
 SDPT3 (fourth column) and AN-FGM (third column) within the the same computational time (fifth column).
The last column indicates the time for AN-FGM to obtain a solution with error up to $0.01\%$
of the final solution generated by SDPT3 (/ indicates that AN-FGM was not able to achieve that accuracy within the allotted time),
that is, ${\|A_{\text{an-fgm}}X-B\|}_F \leq 1.0001 {\|A_{\text{sdpt3}}X-B\|}_F$
where $A_{\text{an-fgm}}$ (resp.\@ $A_{\text{sdpt3}}$) is the solution generated by AN-FGM (resp.\@ SDPT3).

We observe the following:
\begin{itemize}
\item In all the well-conditioned cases and rank-deficient cases, AN-FGM outperforms SDPT3, being order of magnitude faster (comparing the last two columns of Table~\ref{tab:resultsIPM}). This is not surprising since AN-FGM has a much lower per-iteration cost while the convergence will be fast because the problems solved by FGM are well conditioned.
\item  For ill-conditioned cases, SDPT3 allows to obtain high accuracy solutions faster (except in the case $n = 2m = 60$).
However, AN-FGM generates solution at most 2\% from SDPT3 in the worst case.

For $n = 2m = 60$, AN-FGM performs better because the subproblem solved by FGM has size only 30-by-30 hence can perform more iterations.
Moreover, the infimum is not attained which explains why SDPT3 failed to return a solution with acceptable accuracy (\texttt{cvx\_status = Failed} in the 10 cases); see also Remark~\ref{reminfnotat}.
\end{itemize}

\begin{center}
 \begin{table}[h!]
 \begin{center}
\caption{ Comparison between AN-FGM and SDPT3.}
\label{tab:resultsIPM}
 \begin{tabular}{|c|c||c|c|c|c|}
 \hline             &           &   AN-FGM  &  SDPT3 & Time (s.) & AN-FGM - 0.01\% (s.) \\ \hline
 Well-conditioned & $m=n=60$    & 83.77      &  83.77 & 10.75    & 0.36 \\
                  & $m=2n=60$   & 93.06      &  93.06 & 2.26    & 0.01 \\
                  & $n=2m=60$   & 60.85      &  60.86 & 21.38    & 0.02 \\ \hline
 Ill-conditioned  & $m=n=60$    & 78.09      &  \textbf{76.19} & 16.91    & / \\
                  & $m=2n=60$   & 89.44      &  \textbf{88.70} & 2.97    & / \\
                  & $n=2m=60$   & \textbf{54.08}   &  54.36 & 25.93    & 10.78 \\ \hline
 Rank-deficient   & $m=n=60$    & 82.96      &  82.97 & 20.21    & 0.01 \\
                  & $m=2n=60$   & 91.85      &  91.85 & 2.96    & 0.00 \\
                  & $n=2m=60$   & 77.10      &  77.10 & 19.60    & 0.00 \\  \hline
\end{tabular}
 \end{center}
 \end{table}
 \end{center}

\noindent \paragraph{Summary} Except for ill-conditioned problems of relative small size ($m$ and $n$ up to a hundred on our machine) where $n$ is not significantly smaller than $m$, AN-FGM should be preferred to SDPT3.
For large problems, SDPT3 quickly becomes impractical (for example, SDPT3 requires about 80 seconds for the well-conditioned case with $m=n=100$, while Matlab crashed when we tried $m=n=200$).

\begin{remark}{\rm
We also performed a comparison with the interior-point method dedicated to the PSD Procrustes problem
proposed in~\cite{krislock2003numerical} (namely, sdls\_precorr\footnote{Available from \url{https://sites.google.com/site/nathankrislock/}}).
However, it performs in general either similarly as SDPT3 or worse (in particular, for the ill-conditioned case with $n=2m=40$, it returns a solution with relative error higher than 100\%). Therefore, we have not included sdls\_precorr in our comparison.
}
\end{remark}

\section{Conclusion}
In this paper we have first completely described
the set of optimal solutions of the PSD Procrustes problem~\eqref{psdproc} when the infimum is attained, and
a sequence of solutions whose objective function value converges to the infimum when the infimum is not attained.
This description relies on the solution of a smaller PSD Procrustes problem
(where $X$ and $B$ are $r$-by-$r$ matrices with $r = \rank(X)$ and $X$ is diagonal with positive diagonal elements)
whose infimum is always attained.
Then, we have applied an optimal first-order method (namely, the fast gradient method) on the subproblem that is guaranteed to converge linearly with rate $(1-\kappa^{-1})$ where $\kappa = \frac{\sigma_1(X)}{\sigma_r(X)}$.
Moreover, to mitigate the slow convergence in ill-conditioned cases, we proposed a new effective recursive initialization scheme based on a hierarchical decomposition of the problem into well-conditioned subproblems.
Finally, our new method, referred to as AN-FGM, was shown to compete favorably with other first-order methods,
and with a second-order method (namely, SDPT3, a state-of-the-art interior-point method).

\small

\bibliographystyle{siam}
\bibliography{bibliostable}

\end{document}